\documentclass[11pt,a4paper]{article}
\usepackage[a4paper,margin=2.2cm]{geometry}
\usepackage[T1]{fontenc}
\usepackage{amsmath,amssymb,amsthm}
\usepackage{enumitem}
\usepackage[colorlinks=true,linkcolor=blue,citecolor=blue,urlcolor=blue]{hyperref}

\newtheorem{theorem}{Theorem}[section]
\newtheorem{lemma}[theorem]{Lemma}

\newtheorem{corollary}[theorem]{Corollary}
\theoremstyle{definition}
\newtheorem{definition}[theorem]{Definition}
\newtheorem{example}[theorem]{Example}
\newtheorem{question}{Question}
\theoremstyle{remark}
\newtheorem{remark}[theorem]{Remark}

\newcommand{\Q}{\mathcal{Q}}
\newcommand{\K}{\mathcal{K}}
\newcommand{\Pp}{\mathcal{P}}
\newcommand{\V}{\mathcal{V}}
\newcommand{\up}{{\uparrow}}
\newcommand{\dn}{{\downarrow}}
\newcommand{\cl}[1]{\overline{#1}}
\title{\textbf{Retracts of well-filtered co-sober spaces are co-sober}}
\author{Kia Keng Giam\\[2pt]
  \textit{Temasek Junior College, Singapore}\\[2pt]
  Correspondence: \texttt{giamkiakeng@gmail.com}}
\date{September 20, 2026; revised November 9, 2026}

\begin{document}
\maketitle

\begin{abstract}
Xu and Zhao asked whether co-sobriety is inherited by the image of a continuous kernel operator,
and more generally whether it is inherited by retracts. Both questions have recently been answered in
the negative by Xu, Xu, Ji and Qiu. We prove that both questions have positive answers whenever the
ambient space is, in addition, \emph{well-filtered}; in particular they hold for every sober co-sober
space. It follows that the ambient space of any counterexample to either question, such as that of
Xu, Xu, Ji and Qiu, must be co-sober and not well-filtered. We give an example of a kernel-operator
image that is not a direct instance of the inheritance results of Wen and Xu, and we identify the
step of the argument that uses well-filteredness.

\medskip
\noindent\textit{Keywords:} co-sober space; k-irreducible set; well-filtered space; retract; kernel
operator; compact saturated set; sober space.

\medskip
\noindent\textit{2020 MSC}: 54D99; 54B30; 54D30; 06B35.
\end{abstract}

\section{Introduction}

Sobriety is arguably the single most useful property of non-Hausdorff topological spaces, and the
Hofmann--Mislove theorem describes it through an order isomorphism between the compact saturated
sets of a space and the open filters of its lattice of opens. Seeking a dual form of that theorem,
Escard\'o, Lawson and Simpson \cite[Definition~8.9]{ELS04} introduced \emph{k-irreducible} compact
saturated sets and the corresponding class of \emph{co-sober} spaces, and asked whether every sober
space is co-sober \cite[Problem~9.7]{ELS04}. Wen and Xu \cite[Example~2.8]{WX18} answered this in the
negative: Isbell's complete lattice, equipped with the lower topology, is sober but not co-sober.
Shen, Wu, Xi and Zhao \cite{SWXZ20} subsequently produced a dcpo whose Scott space is sober but not
co-sober, and Ji and Xu \cite{JX26} have shown that the Scott space of a complete Boolean algebra need
not be co-sober.

In their survey of open problems on well-filtered and sober spaces, Xu and Zhao \cite[\S9]{XZ21}
collect six questions on co-sobriety. Several have since been settled. The Smyth power space of a
co-sober space need not be co-sober \cite{XXJQ26s}, answering \cite[Question~9.4]{XZ21} negatively;
the category of co-sober spaces is not reflective in the category of $T_0$ spaces \cite{SXZ24},
answering \cite[Question~9.9]{XZ21} negatively; and the two questions concerning inheritance, which
we consider here, have recently been answered negatively by Xu, Xu, Ji and Qiu \cite{XXJQ26k}:

\begin{question}[{\cite[Question~9.5]{XZ21}}]\label{q:kernel}
Let $X$ be co-sober and let $r\colon X\to X$ be a continuous kernel operator for the specialization
order. Is $r(X)$ co-sober?
\end{question}

\begin{question}[{\cite[Question~9.6]{XZ21}}]\label{q:retract}
Is a retract of a co-sober space co-sober?
\end{question}

A positive answer to Question~\ref{q:retract} would give a positive answer to
Question~\ref{q:kernel}. What is known is due to Wen and Xu:
saturated subspaces and closed subspaces of a co-sober space are co-sober
\cite[Theorem~2.9]{WX18}, and so is the image of a continuous \emph{closure} operator (\cite{WX18},
stated as \cite[Theorem~9.2]{XZ21}). The kernel-operator case of Question~\ref{q:kernel} is exactly
the order-dual of that last result, and it is precisely the dual that resists, because co-sobriety
is defined through compact \emph{saturated}, hence upward closed, sets, whereas a kernel operator
moves points downwards. Indeed, Xu, Xu, Ji and Qiu \cite[Theorem~1.1]{XXJQ26k} have constructed a
countably infinite compact co-sober $T_0$ space with a continuous kernel operator whose image is not
co-sober, so both questions have negative answers in general.

Our result is that well-filteredness of the ambient space restores inheritance, already for the
stronger Question~\ref{q:retract}.

\begin{quote}
\textbf{Theorem~\ref{thm:main}.} \textit{Let $X$ be a co-sober, well-filtered $T_0$ space and let
$r\colon X\to X$ be a continuous idempotent. Then the retract $Y=r(X)$ is co-sober.}
\end{quote}

Since every sober space is well-filtered, both questions have positive answers for every sober
co-sober space (Corollary~\ref{cor:sober}), and the ambient space of any counterexample, such as the
one in \cite{XXJQ26k}, must be co-sober and not well-filtered (Corollary~\ref{cor:counterexample}).
Example~\ref{ex:finite} exhibits a kernel-operator image that is neither saturated, nor closed, nor
the image of any closure operator, so it is not a direct instance of the inheritance results of
\cite{WX18}. Remarks~\ref{rem:free} and~\ref{rem:fails} discuss the role of the hypothesis in the
proof: the argument turns on a union-prime up-set of compact saturated sets having a minimal member,
and we illustrate in a co-sober space that such an up-set need not have one.

\noindent\textbf{A note on notation.} The letter $k$ is entrenched in the term ``k-irreducible'' and
is also the letter used for the kernel operator in \cite{XZ21}. To avoid a collision we write $r$ for
the idempotent throughout.

\section{Preliminaries}

All spaces are $T_0$. We follow the conventions of \cite{XZ21}. The \emph{specialization order} of a
space $X$ is
\[
  x\le y \iff x\in\cl{\{y\}},
\]
so that $\cl{\{y\}}=\dn y$, open sets are upper sets, and a subset is \emph{saturated}, that is, an
intersection of open neighbourhoods, exactly when it is an upper set. We write $\Q(X)$ for the set
of nonempty compact saturated subsets of $X$, and $\up_X$, $\dn_X$ for up- and down-closure computed
in $X$. Subspaces carry the subspace topology, and the specialization order of a subspace is the
restriction of the ambient one.

\begin{definition}[Escard\'o--Lawson--Simpson {\cite[Definition~8.9]{ELS04}}]\label{def:kirr}
Let $K\in\Q(X)$. Then $K$ is \emph{k-irreducible} if $K=K_1\cup K_2$ with $K_1,K_2\in\Q(X)$ implies
$K=K_1$ or $K=K_2$. The space $X$ is \emph{co-sober} if every k-irreducible $K\in\Q(X)$ is of the
form $K=\up_X a$ for a (necessarily unique) $a\in X$.
\end{definition}

Recall that $X$ is \emph{well-filtered} if for every filtered family $\K\subseteq\Q(X)$ and every
open $U$, $\bigcap\K\subseteq U$ implies $K\subseteq U$ for some $K\in\K$; see
\cite[Chapter~8]{GL13}. Every sober space is well-filtered, and the converse fails.

A continuous map $r\colon X\to X$ with $r\circ r=r$ is an \emph{idempotent}; its image
$Y=r(X)=\{x\in X:r(x)=x\}$ is precisely a \emph{retract} of $X$. If moreover $r(x)\le x$ for all $x$,
then $r$ is a \emph{kernel operator} for the specialization order; dually, if $x\le r(x)$ for all
$x$, then $r$ is a \emph{closure operator}. For a kernel operator, $r(x)=\max(\dn_X x\cap Y)$; for a
closure operator, $r(x)=\min(\up_X x\cap Y)$.

\begin{lemma}\label{lem:basic}
Let $X$ be a space.
\begin{enumerate}[label=(\roman*),align=right]
  \item Every continuous map is monotone for the specialization order.
  \item If $C\subseteq X$ is compact, then $\up_X C$ is compact saturated.
  \item If $X$ is well-filtered and $\K\subseteq\Q(X)$ is a nonempty filtered family, then
        $\bigcap\K$ is compact saturated.
\end{enumerate}
\end{lemma}

\begin{proof}
(i) is standard. For (ii): saturation is clear. If $\{U_i\}$ is an open cover of $\up_X C$, it
covers $C$, so finitely many $U_{i_1},\dots,U_{i_n}$ cover $C$; being open they are upper sets, so
their union contains $\up_X C$. For (iii): an intersection of upper sets is an upper set. Given an
open cover $\V$ of $\bigcap\K$, put $U=\bigcup\V$. Then $\bigcap\K\subseteq U$, so by
well-filteredness $K\subseteq U$ for some $K\in\K$; as $K$ is compact, finitely many members of $\V$
cover $K$, hence cover $\bigcap\K$.
\end{proof}

\section{The main theorem}

\begin{theorem}\label{thm:main}
Let $X$ be a co-sober, well-filtered $T_0$ space and let $r\colon X\to X$ be a continuous idempotent.
Then the retract $Y=r(X)$, with the subspace topology, is co-sober.
\end{theorem}

\noindent\textbf{Idea of the proof.} Fix a k-irreducible $Q\in\Q(Y)$ and let $K=\up_X Q$ be the
saturation of $Q$ in $X$. Let $\Q_K(X)=\{P\in\Q(X):P\subseteq K\}$, ordered by inclusion. It is
closed under binary unions, so $(\Q_K(X),\cup)$ is a union-semilattice with top element $K$. The
assignment $P\mapsto\up_Y r(P)$ sends members of $\Q_K(X)$ to compact saturated subsets of $Q$ and
preserves finite unions. Hence the family
\[
  \Pp=\{\,P\in\Q_K(X)\;:\;\up_Y r(P)=Q\,\}
\]
of those $P$ that still generate all of $Q$ contains $K$, is an up-set of $\Q_K(X)$, and is
\emph{union-prime}, because $Q$ is k-irreducible:
\[
  P_1\cup P_2\in\Pp\implies P_1\in\Pp\ \text{or}\ P_2\in\Pp\qquad(P_1,P_2\in\Q_K(X)).
\]
(It is in general not an up-set of the whole of $\Q(X)$, since a larger compact saturated set need
not be contained in $K$.) Well-filteredness lets us shrink inside $\Pp$ to a minimal member $P_0$;
union-primeness together with minimality makes $P_0$ k-irreducible, so co-sobriety of $X$ gives
$P_0=\up_X a$; and then $r(a)$ is forced to be the least element of $Q$. The only topological input
beyond compactness is that $r^{-1}(\dn_X q)$ is closed, which is what allows well-filteredness to be
applied.

\begin{proof}
Let $Q\in\Q(Y)$ be k-irreducible. We must produce $y\in Y$ with $Q=\up_Y y$. Put
\[
  K:=\up_X Q.
\]
Since $Q$ is compact in the subspace $Y$ it is compact in $X$, so $K\in\Q(X)$ by
Lemma~\ref{lem:basic}(ii).

\smallskip\noindent\textbf{Step 1: $r(K)=Q$.} Let $x\in K$, say $x\ge q$ with $q\in Q$. By
Lemma~\ref{lem:basic}(i), $r(x)\ge r(q)=q$. Also $r(x)\in Y$, and $Q$ is an upper set of $Y$ because
it is saturated in $Y$; hence $r(x)\in Q$. This gives $r(K)\subseteq Q$, and $Q=r(Q)\subseteq r(K)$
since $r$ fixes $Y$ pointwise.

\smallskip\noindent\textbf{Step 2: the map $\varphi$.} Write $\Q_K(X)=\{P\in\Q(X):P\subseteq K\}$ and
define
\[
  \varphi\colon\Q_K(X)\longrightarrow\Q(Y),\qquad \varphi(P)=\up_Y r(P).
\]
This is well defined: $r(P)$ is a nonempty compact subset of $Y$, so $\up_Y r(P)\in\Q(Y)$ by
Lemma~\ref{lem:basic}(ii) applied in $Y$; and $r(P)\subseteq r(K)=Q$ by Step~1, so
$\varphi(P)\subseteq\up_Y Q=Q$. Moreover $\varphi$ is monotone and
\begin{equation}\label{eq:union}
  \varphi(P_1\cup P_2)=\varphi(P_1)\cup\varphi(P_2),
\end{equation}
since $r(P_1\cup P_2)=r(P_1)\cup r(P_2)$ and $\up_Y$ preserves unions. Finally $\varphi(K)=\up_Y Q=Q$
by Step~1. Set
\[
  \Pp:=\{P\in\Q_K(X):\varphi(P)=Q\},
\]
so that $K\in\Pp$ and in particular $\Pp\neq\emptyset$.

\smallskip\noindent\textbf{Step 3: $\Pp$ is closed under nonempty filtered intersections.} Let
$\{P_j\}_{j\in J}\subseteq\Pp$ be a nonempty filtered family and put $P_\infty=\bigcap_j P_j$. Fix
$q\in Q$ and set
\[
  C_q:=r^{-1}\bigl(\dn_X q\bigr).
\]
Since $\dn_X q=\cl{\{q\}}$ is closed and $r$ is continuous, $C_q$ is closed. Because
$\varphi(P_j)=Q\ni q$, there is $p\in P_j$ with $r(p)\le q$, that is, $P_j\cap C_q\neq\emptyset$, for
every $j$. If we had $P_\infty\cap C_q=\emptyset$, then $P_\infty\subseteq X\setminus C_q$, and
well-filteredness applied to the filtered family $\{P_j\}$ and the open set $X\setminus C_q$ would
yield some $j$ with $P_j\subseteq X\setminus C_q$, contradicting $P_j\cap C_q\neq\emptyset$. Hence
\begin{equation}\label{eq:meet}
  P_\infty\cap C_q\neq\emptyset\qquad\text{for every }q\in Q.
\end{equation}
In particular $P_\infty\neq\emptyset$ (here we use $Q\neq\emptyset$), so $P_\infty\in\Q(X)$ by
Lemma~\ref{lem:basic}(iii); clearly $P_\infty\subseteq K$. By~\eqref{eq:meet}, for each $q\in Q$ there
is $p\in P_\infty$ with $r(p)\le q$, i.e.\ $q\in\up_Y r(P_\infty)=\varphi(P_\infty)$. Thus
$Q\subseteq\varphi(P_\infty)$, and the reverse inclusion holds by Step~2. So $P_\infty\in\Pp$.

\smallskip\noindent\textbf{Step 4: a minimal member.} Order $\Pp$ by $\subseteq$. Every nonempty
chain in $\Pp$ is a nonempty filtered family, so by Step~3 it has a lower bound in $\Pp$, namely its
intersection; the empty chain has the lower bound $K$. By Zorn's Lemma, $\Pp$ has a
$\subseteq$-minimal element $P_0$.

\smallskip\noindent\textbf{Step 5: $P_0$ is k-irreducible in $X$.} Suppose $P_0=P_1\cup P_2$ with
$P_1,P_2\in\Q(X)$. Then $P_i\subseteq P_0\subseteq K$, so $P_i\in\Q_K(X)$ and, by~\eqref{eq:union},
\[
  Q=\varphi(P_0)=\varphi(P_1)\cup\varphi(P_2),
\]
a union of two members of $\Q(Y)$. Since $Q$ is k-irreducible, $\varphi(P_i)=Q$ for some $i$; that
is, $P_i\in\Pp$. As $P_i\subseteq P_0$, minimality of $P_0$ forces $P_i=P_0$.

\smallskip\noindent\textbf{Step 6: conclusion.} $X$ is co-sober and $P_0\in\Q(X)$ is k-irreducible,
so $P_0=\up_X a$ for some $a\in X$; note that $a\in P_0\subseteq K$. If $x\ge a$ then $r(x)\ge r(a)$,
so $r(\up_X a)\subseteq\up_Y r(a)$ and hence
\[
  Q=\varphi(P_0)=\up_Y r(\up_X a)\subseteq\up_Y r(a).
\]
Conversely $r(a)\in r(K)=Q$ by Step~1, and $Q$ is an upper set of $Y$, so $\up_Y r(a)\subseteq Q$.
Therefore $Q=\up_Y r(a)$, which is what we had to prove.
\end{proof}

\section{Consequences}

\begin{corollary}\label{cor:sober}
If $X$ is co-sober and well-filtered, then Questions~\ref{q:kernel} and~\ref{q:retract} have positive
answers for $X$. In particular, every retract of a sober co-sober space is co-sober.
\end{corollary}

\begin{proof}
A continuous kernel operator is in particular a continuous idempotent, so Theorem~\ref{thm:main}
applies; note that the hypothesis $r(x)\le x$ was never used. The last sentence follows because every
sober space is well-filtered.
\end{proof}

\begin{remark}\label{rem:wfretract}
Xu and Zhao proved that every retract of a well-filtered space is well-filtered
\cite[Corollary~3.3]{XZ20}. Combined with Theorem~\ref{thm:main}, this shows that the class of
co-sober well-filtered spaces is closed under retracts.
\end{remark}

\begin{remark}\label{rem:class}
The class covered by Corollary~\ref{cor:sober} is not small. Let $L$ be a complete lattice and let
$\omega(L)$ be the lower topology, generated by the sets $L\setminus\up x$. Then $(L,\omega(L))$ is
sober: indeed $\omega(L)=\nu(L^{\mathrm{op}})$, the upper topology of the dual lattice, and
$L^{\mathrm{op}}$ is again a complete lattice, so every irreducible subset has a supremum and
\cite[Proposition~2.9]{XSXZ20} applies; this is also recorded in \cite{SWXZ20}. By
\cite[Corollary~2.5]{WX18}, $(L,\omega(L))$ is co-sober precisely when the Scott space $\Sigma L$ is
sober. So every complete lattice with sober Scott topology --- for instance every continuous lattice
--- yields a sober co-sober space to which Corollary~\ref{cor:sober} applies. Finite sober co-sober
spaces are exhibited in Example~\ref{ex:finite} below.
\end{remark}

The next example gives a retract to which Theorem~\ref{thm:main} applies but which is not a direct
instance of the three inheritance results of \cite{WX18} recalled in Section~1.

\begin{example}\label{ex:finite}
Let $X=\{\bot,s,c,d\}$ be the poset with $\bot<s<c$ and $s<d$ (so $c$ and $d$ are incomparable),
carrying the Alexandroff topology, whose open sets are the upper sets. Every finite $T_0$ space is
sober, hence well-filtered: an irreducible closed set $A$ equals the finite union
$\bigcup_{a\in A}\dn a$ of closed sets, so irreducibility forces $A=\dn a$ for some $a$. Dually every
finite $T_0$ space is co-sober: a compact saturated $K$ equals the finite union
$\bigcup_{m\in\min K}\up m$ of members of $\Q(X)$, so k-irreducibility forces $K=\up m$ for some $m$.
Thus $X$ satisfies the hypotheses of Theorem~\ref{thm:main}.

Put $Y=\{\bot,c,d\}$ and define $r\colon X\to X$ by $r(\bot)=r(s)=\bot$, $r(c)=c$, $r(d)=d$. Then $r$
is monotone, hence continuous for the Alexandroff topology; it is idempotent; and $r(x)\le x$ for all
$x$, so $r$ is a continuous kernel operator with $r(X)=Y$. Now:
\begin{itemize}
  \item $Y$ is \emph{not saturated}: $\bot\in Y$ and $\bot<s$, but $s\notin Y$.
  \item $Y$ is \emph{not closed}: $c\in Y$ and $s<c$, but $s\notin Y$.
  \item $Y$ is \emph{not the image of any closure operator} on $X$: such an operator would have to
        send $s$ to the least element of $\up_X s\cap Y=\{c,d\}$, and this set has no least element.
\end{itemize}
Hence none of the three inheritance results recalled in Section~1 applies to $Y$, while
Theorem~\ref{thm:main} does. (One checks directly that $Y$, which is the poset $\bot<c$, $\bot<d$ with
the Alexandroff topology, is indeed co-sober.) Infinite examples are obtained in the same way; the
point of taking $X$ finite is that every hypothesis can be verified by inspection.
\end{example}

\begin{corollary}\label{cor:counterexample}
Let $X$ be a co-sober space with a retract that is not co-sober. Then $X$ is \emph{not}
well-filtered. In particular, the co-sober space of \cite[Theorem~1.1]{XXJQ26k}, which has a
continuous kernel operator whose image is not co-sober, is not well-filtered.
\end{corollary}

\begin{proof}
This is the contrapositive of Theorem~\ref{thm:main}.
\end{proof}

\begin{remark}[the hypothesis is not free]\label{rem:free}
Co-sobriety does not imply well-filteredness. For a set $M$ with a topology $\tau$ and $U\in\tau$, let
$\delta_U$ be the topology generated by the cofinite topology on $M$ together with $U$. Then
$\Q((M,\delta_U))=2^M\setminus\{\emptyset\}$, the space $(M,\delta_U)$ is co-sober, and it is
well-filtered if and only if $M$ is finite \cite[Lemma~4.1]{XXJQ26s}. So
Corollary~\ref{cor:counterexample} does not vacuously exclude everything; indeed, the co-sober space of
\cite{XXJQ26k} is not well-filtered and has a retract that is not co-sober.
\end{remark}

\begin{remark}[exactly what fails]\label{rem:fails}
The proof uses co-sobriety only in Step~6 and well-filteredness only in Step~3. What Step~3 buys is
that the union-prime up-set $\Pp$ of the semilattice $(\Q_K(X),\cup)$ has a minimal member. Without
well-filteredness, union-prime up-sets of compact saturated sets need not have minimal members: in
$(M,\delta_U)$ with $M$ infinite, $\Q((M,\delta_U))$ is the semilattice of all nonempty subsets of
$M$, and $\{A\subseteq M:A\text{ infinite}\}$ is a union-prime up-set with no minimal element. This
family only illustrates the minimality issue; it is not shown to arise from a retraction and a
k-irreducible set as in the proof of Theorem~\ref{thm:main}. An actual counterexample to
Questions~\ref{q:kernel} and~\ref{q:retract} is given in \cite{XXJQ26k}, and every counterexample must
exploit precisely this phenomenon: by Steps~5 and~6, which do not use well-filteredness, the family
$\Pp$ attached to a k-irreducible set that is not of the form $\up_Y y$ cannot have a minimal member.
\end{remark}

\begin{remark}\label{rem:notce}
The spaces $(M,\delta_U)$ are themselves not counterexamples: they are $T_1$, so their
specialization order is discrete and they admit no kernel operator other than the identity; and the
subspace induced on any $N\subseteq M$ is again of the form $\delta_{U\cap N}$, hence co-sober.
\end{remark}

\section{Two further observations}

The following is independent of Theorem~\ref{thm:main}, but explains why a non-principal
k-irreducible compact saturated set cannot be downward directed, and so further constrains the
counterexamples to Questions~\ref{q:kernel} and~\ref{q:retract}, such as the one in \cite{XXJQ26k}.

\begin{lemma}\label{lem:directed}
Let $X$ be $T_0$ and let $Q\in\Q(X)$ be downward directed. Then $Q$ has a least element $q_0$, and
$Q=\up_X q_0$.
\end{lemma}

\begin{proof}
Suppose $Q$ has no least element. Given $q\in Q$, choose $p\in Q$ with $p\not\ge q$ and then, by
directedness, $t_q\in Q$ with $t_q\le q$ and $t_q\le p$. If $t_q=q$ then $q\le p$, a contradiction;
so $t_q<q$. As $X$ is $T_0$ and $t_q\le q$ with $t_q\neq q$, there is an open $U_q$ with $q\in U_q$
and $t_q\notin U_q$. Now $\{U_q\}_{q\in Q}$ covers $Q$, so $Q\subseteq U_{q_1}\cup\dots\cup U_{q_n}$
for some $q_1,\dots,q_n$, with $n\ge1$ because $Q\neq\emptyset$. Choose $t\in Q$ below every
$t_{q_i}$. Each $U_{q_i}$ is an upper set omitting $t_{q_i}\ge t$, hence omits $t$. So $t\in Q$ is
uncovered, a contradiction. Finally $Q$ is an upper set containing its least element $q_0$, so
$Q=\up_X q_0$.
\end{proof}

\begin{corollary}\label{cor:total}
Every $T_0$ space whose specialization order is total is co-sober, and so is every subspace of such a
space.
\end{corollary}

\begin{proof}
In a total order every nonempty subset is downward directed, so Lemma~\ref{lem:directed} applies. The
specialization order of a subspace is the restriction of the ambient one, hence again total.
\end{proof}

\subsection*{Acknowledgements}

The author thanks the anonymous referee for a careful reading of the manuscript and for valuable
comments that improved the paper, in particular for pointing out the paper \cite{XZ20} and the
preprint \cite{XXJQ26k}.

\subsection*{CRediT authorship contribution statement}

Kia Keng Giam: Conceptualization, Formal analysis, Investigation, \mbox{Writing -- original draft},
\mbox{Writing -- review \& editing}.

\subsection*{Funding}

This research did not receive any specific grant from funding agencies in the public, commercial, or
not-for-profit sectors.

\subsection*{Declaration of competing interest}

The author declares that they have no known competing financial interests or personal relationships
that could have appeared to influence the work reported in this paper.

\subsection*{Data availability}

No data was used for the research described in the article.

\subsection*{Declaration of generative AI and AI-assisted technologies in the manuscript
preparation process}

During the preparation of this revision the author used Claude (Anthropic), an AI assistant, in order
to check the arguments and references and to draft changes in response to the referee's report. After
using this tool, the author reviewed and edited the content as needed and takes full responsibility
for the content of the published article.

\end{document}